\documentclass[11pt]{article}

\usepackage{amsmath}
\usepackage{amssymb}
\usepackage{latexsym}
\usepackage{comment}
\usepackage{color}
\usepackage{enumitem}
\usepackage{mathrsfs}

\newtheorem{assumption}{Assumption}

\def\qed{ \ \vrule width.2cm height.2cm depth0cm\smallskip}
\newenvironment{proof}{\noindent {\bf Proof.\/}}{$\qed$\vskip 0.1in}

\newcommand{\ba}{\begin{array}}
\newcommand{\ea}{\end{array}}
\newcommand{\be}{\begin{equation}}
\newcommand{\ee}{\end{equation}}
\newcommand{\bea}{\begin{eqnarray}}
\newcommand{\eea}{\end{eqnarray}}
\newcommand{\beaa}{\begin{eqnarray*}}
\newcommand{\eeaa}{\end{eqnarray*}}

\def\dbE{\mathbb{E}}

\def\dbL{\mathbb{L}}

\def\dbP{\mathbb{P}}
\def\dbR{\mathbb{R}}

\def\dbZ{\mathbb{Z}}

\def\d{\delta}
\def\e{\varepsilon}

\def\o{\omega}

\def\O{\Omega}
\def\cF{{\cal F}}

\def\cL{{\cal L}}

\def\cP{{\cal P}}

\def\no{\noindent}

\def\bs{\bigskip}
\def\q{\quad}

\def\qed{ \hfill \vrule width.25cm height.25cm depth0cm\smallskip}

\newcommand{\basa}{\begin{assumption}}
\newcommand{\easa}{\end{assumption}}

\newcommand{\bas}{\begin{assum}}
\newcommand{\eas}{\end{assum}}

\def\limsup{\mathop{\overline{\rm lim}}}
\def\liminf{\mathop{\underline{\rm lim}}}

\def\1{{\bf 1}}

\def\:{\!:\!}
\def\reff#1{{\rm(\ref{#1})}}
\def \proof{{\noindent \bf Proof\quad}}

at 9pt

\begin{document}

\newtheorem{thm}{Theorem}
\newtheorem{lem}[thm]{Lemma}
\newtheorem{cor}[thm]{Corollary}
\newtheorem{prop}[thm]{Proposition}
\newtheorem{rem}[thm]{Remark}
\newtheorem{eg}[thm]{Example}
\newtheorem{defn}[thm]{Definition}
\newtheorem{assum}[thm]{Assumption}

\renewcommand {\theequation}{\arabic{equation}}
\def\thesection{\arabic{section}}

\title{\bf  An Elementary Proof for the Structure of Wasserstein Derivatives}

\author{
Cong Wu \thanks{ \no Department of Mathematics, University of Southern California, Los Angeles, CA 90089. E-mail: congw@usc.edu.} ~ and 
~{Jianfeng Zhang} \thanks{\noindent
Department of Mathematics, University of Southern California, Los Angeles, CA 90089. E-mail: jianfenz@usc.edu. This author is supported in part by NSF grant \#1413717. }
}

\date{}
\maketitle

\begin{abstract}
Let $F: \dbL^2(\O, \dbR)\footnote{The space $\dbR$ can be replaced with general $\dbR^d$. We assume $d=1$ here for simplicity. } \to \dbR$ be a law invariant and continuously Fr\'echet differentiable mapping. Based on Lions \cite{Lions}, Cardaliaguet \cite{Cardaliaguet} (Theorem 6.2 and 6.5) proved that: 
\bea
\label{Derivative}
D F (\xi) = g(\xi),
\eea
where $g: \dbR\to \dbR$ is a deterministic function which depends only on the law  of $\xi$. See also Carmona \& Delarue \cite{CD} Section 5.2 and Gangbo \&  Tudorascu \cite{GT}. In this short note we provide an elementary proof for this well known result.  This note is part of our accompanying paper \cite{WZ}, which deals with a more general situation.
\end{abstract}

\bs

\no Let $\cP_2(\dbR)$ denote the set of square integrable probability measures on $\dbR$, and consider a mapping $f: \cP_2(\dbR) \to \dbR$. As in standard literature, we lift $f$ to  a function $F: \dbL^2(\O, \dbR)\to \dbR$ by $F(\xi):=f(\cL_\xi)$, where  $(\O,\cF,\dbP)$ is an atomless Polish probability space and $\cL_\xi$ denotes the law of $\xi$. 
 If $F$ is Frech\'et differentiable, then $DF(\xi)$ can be identified as an element of $\dbL^2(\O, \dbR)$:
\bea
\label{Frechet}
\dbE\big[ DF(\xi)  \eta\big] = \lim_{\e\to 0} {F(\xi + \e \eta) - F(\xi) \over \e},\q \mbox{for all}~\eta \in \dbL^2(\O, \dbR).
\eea

We start with the simple case that $\xi$ is discrete. Let $\d_x$ denote the Dirac measure of $x$.

\begin{prop}
\label{prop-discrete}
Assume  
$\xi$ is discrete: $\dbP(\xi = x_i) = p_i$, $i\ge 1$.  If $F$ is Fr\'echet differentiable at $\xi$, then \reff{Derivative} holds with
\bea
\label{formula1}
g(x_i) :=  \lim_{\e\to 0} {f(\sum_{j\neq i} p_j \d_{x_j} + p_i \d_{x_i+\e}) - f(\sum_{j\ge 1} p_j \d_{x_j})\over \e p_i}, \q i\ge 1.
\eea
\end{prop}

To prove the proposition, we need the following  result.
\begin{lem}
\label{lem-elem}
Let  
$X\in \dbL^2(\O, \dbR)$.  Assume $A\in \cF$  with $\dbP(A)>0$  satisfies   
\bea
\label{elem}
\dbE[X \1_{A_1}]=\dbE[X \1_{A_2}], \q\mbox{for all}~ A_1,A_2\subset A~ \mbox{such that}~\dbP(A_1) = \dbP(A_2).
\eea
Then $X$ is a constant, $\dbP$-a.s.  in $A$.
\end{lem}
\proof This result is elementary, we nevertheless provide a proof for completeness. 

Assume the result is not true. Denote $c := {\dbE[X\1_A]\over \dbP(A)}$ and $A_1 := \{X <c\}\cap A$, $A_2:= \{X>c\}\cap A$.  Then $\dbP(A_1)>0$, $\dbP(A_2) >0$. Assume without loss of generality that $\dbP(A_1) \le \dbP(A_2)$. Since $(\O, \cF, \dbP)$ is atomless,  there is a random variable $U$ with uniform distribution  on $[0, 1]$. Denote $A_{2,x} := A_2 \cap \{U \le x\}$, $x\in [0, 1]$. Clearly there exists $x_0$ such that $\dbP(A_{2,x_0}) = \dbP(A_1)$.  Apply \reff{elem} on $A_1$ and $A_{2,x_0}$ we obtain the desired contradiction.
\qed

\begin{rem}
\label{rem-elem}
{\rm Lemma \ref{lem-elem} may not hold if $(\O, \cF, \dbP)$ has atoms. Indeed, consider $\O := \{\o_1, \o_2\}$ with $\dbP(\o_1) = {1\over 3}, \dbP(\o_2) = {2\over 3}$. Set $A:= \O$ and $X$ is an arbitrary random variable. The \reff{elem} holds true trivially because $\dbP(A_1) \neq \dbP(A_2)$ whenever $A_1\neq A_2$. However, $X$ may not be a constant.
\qed}
\end{rem}

\no{\bf Proof of Proposition \ref{prop-discrete}.} Fix an $i\ge 1$. For an arbitrary $A_1 \subset A:= \{\xi = x_i\}$, set $\eta := \1_{A_1}$.  Note that, for any $\e > 0$,  we have 
\beaa
\cL_{\xi + \e \eta} = \sum_{j\neq i} p_j \d_{x_j} + \dbP(A_1) \d_{x_i + \e} + [p_i - \dbP(A_1)] \d_{x_i},
\eeaa
which depends only on $\cL_\xi$ and $\dbP(A_1)$. By \reff{Frechet},
\bea
\label{formula0}
\dbE\big[DF(\xi) \1_{A_1}\big] =  \lim_{\e\to 0} {f\big( \sum_{j\neq i} p_j \d_{x_j} + \dbP(A_1) \d_{x_i + \e} + [p_i - \dbP(A_1)] \d_{x_i}\big) - f(\sum_{j\ge 1} p_j \d_{x_j})\over \e}.
\eea
In particular, $\dbE\big[DF(\xi) \1_{A_1}\big]$ depends only on $\dbP(A_1)$ for $A_1\subset  \{\xi = x_i\}$. Applying Lemma  \ref{lem-elem}, we see that $DF(\xi)$ is a constant, $\dbP$-a.s. on $\{\xi=x_i\}$. Now set $A_1 := \{\xi=x_i\}$ in \reff{formula0}, we obtain \reff{formula1} immediately. 
\qed

We now consider the general case.

\begin{thm}
\label{thm-general}
If $F$ is continuously Fr\'echet differentiable, then \reff{Derivative} holds with $g$ depending only on $\cL_\xi$ but not on the particular choice of $\xi$. 
\end{thm}
\proof For each $n\ge 1$, denote $x^n_i := i 2^{-n}, i \in \dbZ$, and $\xi_n :=  \sum_{i=-\infty}^\infty x^n_i \1_{\{x^n_i \le \xi < x^n_{i+1}\}}$. Since $\xi_n$ is discrete, by Proposition \ref{prop-discrete} we have $DF(\xi_n) = g_n(\xi_n)= \tilde g_n(\xi)$, where $g_n$ is defined on $\{x^n_i, i\in \dbZ\}$ by \reff{formula1} (with $g_n(x^n_i) :=0$ when $\dbP(\xi_n = x^n_i)=0$) and $\tilde g_n(x) := g_n(x^n_i)$ for $x\in [x^n_i, x^n_{i+1})$.   Clearly $\lim_{n\to \infty}\dbE[|\xi_n-\xi|^2] = 0$. Then by the continuous differentiability of $F$ we see that $\lim_{n\to\infty} \dbE[|\tilde g_n(\xi) - DF(\xi)|^2] = 0$.  Thus, there exists a subsequence $\{n_k\}_{k\ge 1}$ such that $\tilde g_{n_k} (\xi) \to D F(\xi)$, $\dbP$-a.s. Denote $K := \{x: \limsup_{k\to\infty} \tilde g_{n_k}(x) = \liminf_{k\to\infty} \tilde g_{n_k}(x)\}$, and  $g(x) := \lim_{k\to \infty} \tilde g_{n_k}(x)\1_K(x)$. Then $\dbP(\xi\in K)=1$ and  $D F(\xi) = g(\xi)$, $\dbP$-a.s.  

Moreover, let $\xi'$ be another random variable such that $\cL_{\xi'}=\cL_\xi$. Define $\xi_n'$ similarly. Then $D F(\xi'_n) = \tilde g_n(\xi')$ for the same function $\tilde g_n$. Note that $\dbP(\xi'\in K) = \dbP(\xi\in K) =1$, then $\lim_{k\to\infty} \tilde g_{n_k}(\xi') = g(\xi')$, $\dbP$-a.s. On the other hand, $D F(\xi'_{n_k}) \to DF (\xi')$ in $\dbL^2$. So $D F(\xi') = g(\xi')$, and thus $g$ does not depend on the choice of $\xi$.
\qed

\begin{rem}
\label{rem-joint}
{\rm One may also write $D F(\xi) = g(\cL_\xi, \xi)$, where $g: \cP_2(\dbR) \times \dbR \to \dbR$. When $DF$ is uniformly continuous, one may easily construct $g$ jointly measurable in $(\mu, x) \in  \cP_2(\dbR) \times \dbR$. One may also extend the result to the case that $F$ is a function of processes.  We leave the details  to \cite{WZ}.
\qed}
\end{rem}


\begin{thebibliography}{1}

\bibitem{Cardaliaguet}
Cardaliaguet P. (2013), {\it Notes on Mean Field Games (from P.-L. LionsÕ lectures at College de France)}, preprint, www.ceremade.dauphine.fr/$\sim$cardalia/MFG100629.pdf.

\bibitem{CD}
Carmona, R. and  Delarue, F. (2017) {\sl Probabilistic Theory of Mean Field Games I Ð Mean Field FBSDEs, Control, and Games}, Springer Verlag, 2017.

\bibitem{GT}
Gangbo, W.  and  Tudorascu, A. {\it On differentiability in the Wasserstein space and well-posedness for Hamilton-Jacobi equations}. Technical report, 2017.


\bibitem{Lions}
 Lions, P.-L. {\sl Cours au Coll\'ege de France}. www.college-de-france.fr.
 
 
 \bibitem{WZ}
Wu, C. and Zhang, J. {\it Viscosity Solutions to Parabolic Master Equations and McKean-Vlasov SDEs with Closed-loop Controls}, preprint, arXiv:1805.02639.
 
 \end{thebibliography}
\end{document}